# Interval-based parameter identification for structural static problems

N. Xiao[1], F. Fedele[1,2], and R.L. Muhanna[1]

[1] *School of Civil and Environmental Engineering,*
*Georgia Institute of Technology, Atlanta, GA 30332, USA*

[2] *School of Electrical and Computer Engineering,*
*Georgia Institute of Technology, Atlanta, GA 30332, USA*

**Abstract**

We present an interval-based approach for parameter identification in structural static problems. The proposed inverse formulation models uncertainties in measurement data as intervals, and exploits the Interval Finite Element Method (IFEM) combined with adjoint-based optimization. The inversion consists of a two-step algorithm: first, an estimate of the parameters is obtained by means of a deterministic iterative solver. Then, the algorithm switches to the interval extension of the previous solver, using the deterministic estimate of the parameters as an initial guess. Various numerical examples show that the proposed method provides guaranteed interval enclosures of the parameters, and it always contains Monte Carlo predictions.

## 1    Introduction

Parameter identification aims at estimating modal parameters of a physical system from available measurements of the system response. It belongs to the class of inverse problems (e.g. Hansen, 2010; Ramm, 2005; Santamarina and Fratta, 2005). For example, wave tomography is used in geophysics for seismic waveform inversion (Fichtner, 2010); in biomedical engineering, optical tomography is used to detect breast cancer tissue via fluorescence (Fedele, *et al*. 2002; Eppstein, *et al*. 2003); in civil engineering, inversion techniques are used for structural health monitoring or damage detection in safety evaluation (Chang, *et al*. 2003; Glaser, *et al*. 2007). In the abovementioned problems, the system response is predicted based on initial guessed modal parameters, and it is then compared with the actual measurement data. Then, iterative corrections of the modal parameters lead to a solution, which minimizes the difference between the predicted system response and measurement data in a least-square sense.

Inevitably, data contain errors caused by measurement devices or unfriendly environmental conditions during data acquisition. Such uncertainties can be modeled using probability theory (e.g. Akashi, *et al*. 1979; Stull, *et al*. 2011; Wang, *et al*. 2011). For example, Kalman filtering (see Kalman, 1960; Brown and Hwang, 1992; Simon, 2006) provides error estimates on the modal parameters based on noisy measurements of the response of a time-evolving system (e.g. Xie and Feng, 2012). Clearly, probability approaches have their limitations, since they require a prior assumption on the nature of the uncertainty, which is usually modeled as a random Gaussian variable. However, such an assumption is too optimistic or not realistic. In practice, there are often not enough measurements to reliably assess the statistical nature of the associated uncertainties. Instead, we only know bounds on the uncertain variable and some partial information about its probabilities. In this setting, non-probability theories such as fuzzy sets (Haag, *et al*. 2010; Erdogan and Bakir, 2013), evidence theory (Jiang, *et al*. 2013) and intervals (Khodaparast, *et al*. 2011; Du, *et al*. 2014) are useful for modeling uncertainties.

In this work, we exploit the Interval Finite Element Method (IFEM) (Rama Rao, *et al*. 2011; Muhanna *et al*. 2007) combined with adjoint-based optimization (Fedele, *et al*. 2002; Eppstein, *et al*. 2003) to provide a new algorithm that guarantees interval enclosure of the modal parameters from inversion of noisy measurements modeled as intervals.

The paper is organized as follows. First, IFEM is briefly reviewed and new decomposition strategies are presented to limit overestimation due to multiple occurrences of the same variable in the IFEM matrix equations. Then, the deterministic inverse algorithm is formulated using adjoint-based methods. The extension of the algorithm to intervals is then presented. Finally, several numerical examples are discussed to validate the performance of the proposed method.

## 2 Interval Finite Element Method

Interval Finite Element Method (IFEM) uses intervals to describe uncertain variables and follow the general procedure of conventional Finite Element Method (FEM). Intervals are extension of real numbers. Instead of representing one single point in the real axis, an interval denotes set of real numbers, which are most suitably described by its endpoints,

$$\mathbf{x} = [\underline{x}, \bar{x}] = \{x \mid \underline{x} \leq x \leq \bar{x}, x \in \Re\}, \tag{1.}$$

where $\mathbf{x}$ denotes the interval, $\underline{x}$ and $\bar{x}$ denote its lower and upper bounds, respectively, and bold symbols denote interval quantities. Alternatively, an interval can be represented by its midpoint $x_{\mathrm{mid}} = (\underline{x} + \bar{x})/2$ and radius $x_{\mathrm{rad}} = (\bar{x} - \underline{x})/2 \geq 0$. The width of an interval is defined as $x_{\mathrm{wid}} = (\bar{x} - \underline{x}) = 2x_{\mathrm{rad}}$. Intervals with non-zero midpoint values can be brought into the form of $\mathbf{x} = x_{\mathrm{mid}}(1 + \boldsymbol{\delta}_x)$, where $\boldsymbol{\delta}_x$ has a zero midpoint. The width of $\boldsymbol{\delta}_x$ in percentage is usually referred to as the uncertainty level of $\mathbf{x}$. For a detailed discussion on interval arithmetic and extensions to interval matrices and functions, we refer to (Alefeld and Herzberger, 1983; Kulisch and Miranker, 1981; Moore *et al*. 2009).

Overestimation due to dependency is the curse in any implementation of interval arithmetic (see Muhanna and Zhang, 2007; Muhanna, *et al*. 2007). In order to reduce it, we propose a new decomposition strategy for the stiffness matrix $\mathbf{K}$ and the nodal equivalent load $\mathbf{f}$ of a structural system governed by the equilibrium condition $\mathbf{Ku} = \mathbf{f}$. Here, $\mathbf{K}$ and $\mathbf{f}$ are decomposed as

$$\mathbf{K} = A\,\mathrm{diag}(\Lambda\boldsymbol{\alpha})A^T, \quad \mathbf{f} = M\boldsymbol{\delta}, \tag{2.}$$

where $A$, $\Lambda$ and $M$ are scalar matrices; $\boldsymbol{\alpha}$ and $\boldsymbol{\delta}$ are interval vectors containing all the uncertainties in the system; and diag($\mathbf{v}$) maps a vector $\mathbf{v}$ into a diagonal matrix, whose diagonal is $\mathbf{v}$. In this way, we separate deterministic and uncertain terms, and multiple occurrences of the same variable are avoided. In practice, the decomposition in Eq. (2.) is done in two steps. In the first step, the element stiffness matrix $\mathbf{K}_e$ and the element nodal equivalent load $\mathbf{f}_e$ are decomposed into $A_e$, $\Lambda_e$, $M_e$, $\boldsymbol{\alpha}_e$ and $\boldsymbol{\delta}_e$ using Eq. (2) in the local reference system. In the second step, $A_e$, $\Lambda_e$ and $M_e$ are assembled into $A$, $\Lambda$ and $M$ in the global reference system.

In particular, for an element with uncertain material properties

$$\mathbf{K}_e = \int_\Omega B_e^T \mathbf{E}_e B_e d\Omega, \qquad (3.)$$

where the integration domain $\Omega$ is the entire element, $B_e$ is the scalar strain-displacement matrix at arbitrary locations inside the element, and $\mathbf{E}_e$ is the interval constitutive matrix, which is function of material uncertainties. To reduce overestimation due to dependency, $\mathbf{K}_e$ is decomposed as

$$\mathbf{K}_e = A_e \mathrm{diag}(\Lambda_e \boldsymbol{\alpha}_e) A_e^T, \qquad (4.)$$

where $A_e$, $\Lambda_e$ are scalar matrices, and the interval vector $\boldsymbol{\alpha}_e$ contains all the uncertainties of the element.

From Eq. (3), numerical integration yields

$$\mathbf{K}_e = \sum_{i=1}^{m} w_i J(\xi_i) B_e^T(\xi_i) \mathbf{E}_e(\xi_i) B_e(\xi_i) \qquad (5.)$$

where $m$ is the number of integration points used, $\xi_i$ and $w_i$ are respectively the coordinates and weights of the integration points, and $J$ is the determinant of the Jacobian. The scalar matrices $A_e$, $\Lambda_e$ and the interval vector $\boldsymbol{\alpha}_e$ in Eq. (4.) are given by

$$A_e = \left\{ B_e^T(\xi_1)\Phi_e \quad \cdots \quad B_e^T(\xi_m)\Phi_e \right\},$$
$$\Lambda_e = \begin{Bmatrix} w_1 J(\xi_1)\varphi_e & & \\ & \ddots & \\ & & w_m J(\xi_m)\varphi_e \end{Bmatrix}, \quad \boldsymbol{\alpha}_e = \begin{Bmatrix} \mathbf{E}(\xi_1) \\ \vdots \\ \mathbf{E}(\xi_m) \end{Bmatrix}, \qquad (6.)$$

where $\mathbf{E}(\xi_i)$ denote interval Young's modulus at the $i$-th integration point. Further, $\Phi_e$ and $\varphi_e$ are obtained from the interval constitutive matrix, which is decomposed as

$$\mathbf{E}_e(\xi_i) = \Phi_e \mathrm{diag}\{\varphi_e \mathbf{E}(\xi_i)\} \Phi_e^T. \qquad (7.)$$

The decomposition of the element nodal equivalent load $\mathbf{f}_e$ is done exploiting the $M$-$\boldsymbol{\delta}$ method (Mullen and Muhanna, 1999). Here, $\mathbf{f}_e = M_e \boldsymbol{\delta}_e$, or equivalently

$$\mathbf{f}_e = \sum_{i=1}^{n} N^T(\xi_i)\mathbf{f}_c(\xi_i) + \int_\Omega N^T(\xi)\mathbf{f}_d(\xi) d\Omega. \qquad (8.)$$

where $n$ is the number of concentrated loads acting on the element, $N(\xi)$ is the displacement interpolation matrix for the element, $\mathbf{f}_c(\xi_i)$ are the concentrated loads under consideration, $\Omega$ is the integration domain in which the distributed load $\mathbf{f}_d(\xi)$ is non-zero.

A further simplification can be obtained by rewriting $\mathbf{f}_c(\xi) = L_c(\xi)\boldsymbol{\delta}_e$ and $\mathbf{f}_d(\xi) = L_d(\xi)\boldsymbol{\delta}_e$ as function of the interval vector $\boldsymbol{\delta}_e$, where $L_c(\xi)$ and $L_d(\xi)$ are scalar matrices. Then from Eq. (8)

$$\mathbf{f}_e = \sum_{i=1}^{n} N^T(\xi_i) L_c(\xi_i)\boldsymbol{\delta}_e + \int_{\Omega} N^T(\xi) L_c(\xi)\boldsymbol{\delta}_e d\Omega$$
$$= \left\{\sum_{i=1}^{n} N^T(\xi_i) L_c(\xi_i) + \int_{\Omega} N^T(\xi) L_c(\xi) d\Omega\right\}\boldsymbol{\delta}_e = M_e \boldsymbol{\delta}_e. \quad (9.)$$

Here, $M_e$ is the matrix within braces, which depends on the displacement interpolation matrix $N(\xi)$ and load distribution functions $L_c(\xi)$ and $L_d(\xi)$.

The global $\mathbf{K}$ and $\mathbf{f}$ follow from the conventional assembly strategy (Cook, *et al*. 2007), i.e.

$$\mathbf{K} = \sum_e T_e^T \mathbf{K}_e T_e, \quad \mathbf{f} = \sum_e T_e^T \mathbf{f}_e. \quad (10.)$$

where $T_e$ is the matrix in the transformation $\mathbf{u}_e = T_e \mathbf{u}$ between the global and local nodal displacement vector $\mathbf{u}$ and $\mathbf{u}_e$. Note that $\mathbf{K}_e$, $\mathbf{f}_e$ and $T_e$ are not necessarily the same for each element. By inserting $\mathbf{K}_e = A_e \mathrm{diag}(\Lambda_e \boldsymbol{\alpha}_e) A_e^T$ of Eq. (4.) into Eq. (10.), the decomposition rule for $\mathbf{K}$ follows as

$$\mathbf{K} = \sum_e T_e^T \mathbf{K}_e T_e = \sum_e T_e^T A_e \mathrm{diag}(\Lambda_e \boldsymbol{\alpha}_e) A_e^T T_e$$
$$= \left\{T_e^T A_e \quad \cdots \quad T_e^T A_e\right\} \begin{bmatrix} \mathrm{diag}(\Lambda_e \boldsymbol{\alpha}_e) & & \\ & \ddots & \\ & & \mathrm{diag}(\Lambda_e \boldsymbol{\alpha}_e) \end{bmatrix} \begin{Bmatrix} A_e^T T_e \\ \vdots \\ A_e^T T_e \end{Bmatrix}. \quad (11.)$$

Here, the vector $\boldsymbol{\alpha}_e$ lists the uncertain interpolated Young's modulus at the element integration points, and it relates to the system parameter vector $\boldsymbol{\alpha}$ via $\boldsymbol{\alpha}_e = L_\alpha \boldsymbol{\alpha}$. Comparing terms in Eqs. (2.) and (11.) yields the assembly rules for $A$ and $\Lambda$

$$A = \left\{T_e^T A_e \quad \cdots \quad T_e^T A_e\right\}, \quad \Lambda = \begin{Bmatrix} \Lambda_e & & \\ & \ddots & \\ & & \Lambda_e \end{Bmatrix} \begin{Bmatrix} L_\alpha \\ \vdots \\ L_\alpha \end{Bmatrix} = \begin{Bmatrix} \Lambda_e L_\alpha \\ \vdots \\ \Lambda_e L_\alpha \end{Bmatrix}. \quad (12.)$$

Again, note that $A_e$, $\Lambda_e$ and $L_\alpha$ are not necessarily the same for each element. Similarly, the decomposition rule for $\mathbf{f}$ and the assembly rule for $M$ follow by introducing $\mathbf{f}_e = M_e \boldsymbol{\delta}_e$ into Eq. (10.) and setting $\boldsymbol{\delta}_e = L_\delta \boldsymbol{\delta}$, that is

$$\mathbf{f} = \sum_e T_e^T \mathbf{f}_e = \sum_e T_e^T M_e \boldsymbol{\delta}_e = \left\{T_e^T M_e \quad \cdots \quad T_e^T M_e\right\} \begin{Bmatrix} \boldsymbol{\delta}_e \\ \vdots \\ \boldsymbol{\delta}_e \end{Bmatrix};$$
$$\Rightarrow M = \left\{T_e^T M_e \quad \cdots \quad T_e^T M_e\right\} \begin{Bmatrix} L_\delta \\ \vdots \\ L_\delta \end{Bmatrix} = \sum_e T_e^T M_e L_\delta. \quad (13.)$$

The resulting stiffness matrix **K** in Eq. (10.) is still singular, as essential boundary conditions have not been applied yet. To eliminate the singularity, **u** must satisfy the additional constraint $C\mathbf{u} = 0$, with $C$ denoting a constraint matrix (Rama Rao, *et al*. 2011). Each row of $C$ states one constraint, and the corresponding entry is set equal to 1, leaving the rest of the row null. Then the equilibrium equation follows from setting to zero the first variation of the energy functional **Π** of the structure

$$\mathbf{\Pi} = \frac{1}{2}\mathbf{u}^T\mathbf{K}\mathbf{u} - \mathbf{u}^T\mathbf{f} + \lambda^T(C\mathbf{u}), \tag{14.}$$

that is

$$\begin{bmatrix} \mathbf{K} & C^T \\ C & 0 \end{bmatrix} \begin{Bmatrix} \mathbf{u} \\ \lambda \end{Bmatrix} = \begin{Bmatrix} \mathbf{f} \\ \mathbf{0} \end{Bmatrix}, \tag{15.}$$

where the Lagrangian multiplier $\lambda$ enforces $C\mathbf{u} = 0$. If **K** is composed of degenerated intervals (intervals with zero width), we can establish a direct relationship between **u** and **f** by inverting the generalized stiffness matrix in Eq. (15.), that is

$$\begin{bmatrix} \mathbf{K} & C^T \\ C & 0 \end{bmatrix}^{-1} = \begin{Bmatrix} \mathbf{G}_{11} & \mathbf{G}_{12}^T \\ \mathbf{G}_{12} & \mathbf{G}_{22} \end{Bmatrix}, \quad \Rightarrow \quad \mathbf{u} = \mathbf{G}_1\mathbf{f} = \mathbf{K}^{-1}\mathbf{f}. \tag{16.}$$

In other words, we find the inverse of **K** under the constraint $C\mathbf{u} = 0$ and $\mathbf{K}^{-1} = \mathbf{G}_{11}$.

## 3 Deterministic Inverse Solver

Given an interval load uncertainty vector **δ** and a measurement vector **ũ**, a deterministic solution of the modal parameters **α** is sought using midpoint values of **δ** and **ũ**, and all interval quantities are replaced with their midpoint values. Drawing from Fedele, *et al*. (2014), the algorithm is derived using adjoint based optimization and it exploits conjugate gradient type methods to find optimal estimates of the unknown parameters.

Assume measurements *ũ* are collected at sampling points on the structure and given in terms of the nodal displacement vector *u*, viz. *ũ* = *Hu*. The proposed inverse solver aims at minimizing the difference between the predicted response *Hu* and the actual measurements vector *ũ*, under the equilibrium constraint *Ku* = *f*. To do so, define the objective functional

$$\Gamma = \frac{1}{2}(Hu-\tilde{u})^T S(Hu-\tilde{u}) + w^T(Ku-f) + \frac{1}{2}\alpha^T(\gamma R)\alpha, \tag{17.}$$

where *S* is a diagonal matrix defining the weight for each measurement, *w* is the Lagrangian multiplier to enforce equilibrium (Fedele, *et al.* 2014) and the last term provides regularization for the problem if necessary (e.g. Hansen, 2010). Here, *γ* is the regularizer weight and *R* is the finite-difference matrix associated with second-order differentiation (e.g. Hansen, 2010; Santamarina and Fratta, 2005).

From the decomposition in Eqs. (2) and (17), the first variation of Γ

$$\delta\Gamma = \delta u^T H^T S(Hu - \tilde{u}) + \delta w^T (Ku - f) + \delta u^T K^T w$$
$$+ w^T A \text{diag}(\Lambda \delta\alpha) A^T u + \delta\alpha^T (\gamma R)\alpha. \quad (18.)$$

is null if

$$\begin{cases} Ku - f = 0; \\ Kw + (H^T SH)u - (H^T S)\tilde{u} = 0; \\ \Lambda^T (A^T w \circ A^T u) + (\gamma R)\alpha = 0, \end{cases} \quad (19.)$$

where $a \circ b$ denotes the element-by-element (Hadamard) product of two vectors $a$ and $b$. To obtain Eq. (19) from Eq. (18), we have exploited the matrix symmetry [see Eq. (2)]

$$K = K^T = A\text{diag}(\Lambda\alpha)A^T, \quad (20.)$$

and the chain of identities

$$w^T A\text{diag}(\Lambda\delta\alpha)A^T u = w^T A(\Lambda\delta\alpha \circ A^T u) = \delta\alpha^T \Lambda^T (A^T w \circ A^T u) \quad (21.)$$

The three equations in Eq. (19.) can be interpreted as: *i*) equilibrium condition of the original system with equivalent load *f*, *ii*) equilibrium condition for the adjoint system with equivalent load $B^T S(\tilde{u} - Bu)$, and *iii*) optimal condition that the gradient *g* of $\Gamma$ with respect to $\alpha$ is zero at the solution point.

The first two equations in Eq. (19), viz. the equilibrium conditions for the original and adjoint systems, can be recast in block form

$$\begin{bmatrix} H^T SH & K \\ K & 0 \end{bmatrix} \begin{Bmatrix} u \\ w \end{Bmatrix} = \begin{bmatrix} 0 & H^T S \\ M & 0 \end{bmatrix} \begin{Bmatrix} \delta \\ \tilde{u} \end{Bmatrix}, \quad (22.)$$

where the decomposition $f = M\delta$ is used. The unknown vectors *u* and *w* follow as

$$\begin{Bmatrix} u \\ w \end{Bmatrix} = \begin{bmatrix} 0 & K^{-1} \\ K^{-1} & -K^{-1} H^T SHK^{-1} \end{bmatrix} \begin{bmatrix} 0 & H^T S \\ M & 0 \end{bmatrix} \begin{Bmatrix} \delta \\ \tilde{u} \end{Bmatrix}. \quad (23.)$$

The corresponding objective functional $\Gamma$ and its gradient *g* with respect to $\alpha$, viz. third equation in Eq. (19), can be expressed in terms of *u*, *w* and $\alpha$ respectively as

$$\Gamma = \frac{1}{2}(Hu - \tilde{u})^T S(Hu - \tilde{u}) + \frac{1}{2}\alpha^T (\gamma R)\alpha;$$

$$g = \frac{\partial \Gamma}{\partial \alpha} = \Lambda^T (A^T w \circ A^T u) + (\gamma R)\alpha. \quad (24.)$$

The conjugate gradient method (Andrei, 2009; Yu, *et al*. 2009; Zhang and Li, 2011) is then exploited to iteratively solve for Eq. (19.). We start from a random initial guess $\alpha_1$ and a descending direction $d_1$ along which $\Gamma$ decreases. A natural choice for $d_1$ is the opposite gradient direction, $d_1 = -g_1$. At the *i*-th step, the modal parameter $\alpha$ is updated as

$$\alpha_{i+1} = \alpha_i + s_i d_i, \tag{25.}$$

where $s_i$ is the step size. Here the inexact line search method is used to find an acceptable $s_i$ along the descending direction $d_i$. This should be large enough to yield a significant decrease in $\Gamma$, while not too large to deviate too far from the optimal point. In the proposed method, we adopt the weak Wolfe criterion (Shi and Shen, 2004; Han, *et al.* 2010)

$$\tau_l \leq \frac{\Gamma(\alpha_{i+1}) - \Gamma(\alpha_i)}{s_i g_i^T d_i}, \quad \frac{g(\alpha_{i+1})^T d_i}{g_i^T d_i} \leq \tau_u \tag{26.}$$

where $0 < \tau_l < \tau_u < 1$. In the next iteration step, the descending direction $d_{i+1}$ is determined by the following iterative rule

$$d_{i+1} = -g_{i+1} + \theta_i d_i, \tag{27.}$$

where the parameter $\theta_i$ can be chosen in various ways. Popular choices for $\theta_i$ include (Hestenes and Stiefel, 1952; Fletcher and Reeves, 1964; Polak and Ribière, 1969; Polyak, 1969)

$$\begin{aligned}
\theta_i &= \frac{g_{i+1}^T (g_{i+1} - g_i)}{d_i^T g_i}, & \text{Hestenes-Stiefel} \\
\theta_i &= \frac{g_{i+1}^T g_{i+1}}{g_i^T g_i}, & \text{Fletche-Revees} \\
\theta_i &= \frac{g_{i+1}^T (g_{i+1} - g_i)}{g_i^T g_i}, & \text{Polak-Ribière-Polyak}
\end{aligned} \tag{28.}$$

The algorithm stops when the gradient $g$ and the update on $\alpha$ are both small enough, that is

$$\|\alpha_{i+1} - \alpha_i\| / \|\alpha_i\| \leq \tau, \quad \|g_{i+1}\| / \|g_1\| \leq \tau, \tag{29.}$$

where $\tau$ is the error tolerance. Hereafter, we will adopt the Polak-Ribière-Polyak rule.

## 4  Interval Inverse Solver

The interval algorithm consists of two steps. In the first step, deterministic solutions $u_0$, $w_0$ and $\alpha_0$ are obtained using the deterministic inverse solver described in the previous section. In the second step, these solutions are used as initial guesses for an interval-based inverse solver, generalization of the deterministic one to intervals. This is formulated drawing from Fedele *et al.* (2014). In particular, given an interval load uncertainty vector $\tilde{\boldsymbol{\delta}}$ and measurements $\tilde{\mathbf{u}}$, the unknown interval $\mathbf{u}$, $\mathbf{w}$ and $\boldsymbol{\alpha}$ satisfy the interval extension of Eq. (19), that is

$$\begin{cases} \mathbf{K}(\boldsymbol{\alpha})\mathbf{u} - M\boldsymbol{\delta} = 0; \\ \mathbf{K}(\boldsymbol{\alpha})\mathbf{w} + (H^T S H)\mathbf{u} - (H^T S)\tilde{\mathbf{u}} = 0; \\ \Lambda^T (A^T \mathbf{w} \circ A^T \mathbf{u}) + (\gamma R)\boldsymbol{\alpha} = 0, \end{cases} \tag{30.}$$

where $\mathbf{K}(\boldsymbol{\alpha})$ emphasizes the dependence on the unknown parameter $\boldsymbol{\alpha}$. In order to solve for Eq. (30.), define $\delta_0$ and $\tilde{u}_0$ as the midpoint values of $\boldsymbol{\delta}$ and $\tilde{\mathbf{u}}$, respectively. Then $\delta_0$, $\tilde{u}_0$, $u_0$, $w_0$ and $\alpha_0$ satisfy the optimality conditions in Eq. (19.). Now, introduce the auxiliary variables

$$\Delta\boldsymbol{\delta} = \delta_0 - \boldsymbol{\delta}, \quad \Delta\tilde{\mathbf{u}} = \tilde{u}_0 - \tilde{\mathbf{u}}, \quad \Delta\mathbf{u} = u_0 - \mathbf{u}, \quad \Delta\mathbf{w} = w_0 - \mathbf{w}, \quad \Delta\boldsymbol{\alpha} = \alpha_0 - \boldsymbol{\alpha} \qquad (31.)$$

to represent deviations of the reference solutions from the corresponding interval vectors. Then, the following equalities hold

$$\begin{aligned}
\mathbf{K}\mathbf{u} &= K_0 u_0 - K_0 \Delta\mathbf{u} - \Delta\mathbf{K} u_0 + \Delta\mathbf{K}\Delta\mathbf{u}; \\
\mathbf{K}\mathbf{w} &= K_0 w_0 - K_0 \Delta\mathbf{w} - \Delta\mathbf{K} w_0 + \Delta\mathbf{K}\Delta\mathbf{w}; \\
A^T \mathbf{w} \circ A^T \mathbf{u} &= A^T w_0 \circ A^T u_0 - A^T w_0 \circ A^T \Delta\mathbf{u} \\
&\quad - A^T \Delta\mathbf{w} \circ A^T u_0 + A^T \Delta\mathbf{w} \circ A^T \Delta\mathbf{u}.
\end{aligned} \qquad (32.)$$

These together with

$$\Delta\mathbf{K}\mathbf{u} = A\,\mathrm{diag}(\Lambda\Delta\boldsymbol{\alpha})A^T \mathbf{u} = A\,\mathrm{diag}(A^T\mathbf{u})\Lambda\Delta\boldsymbol{\alpha} \qquad (33.)$$

are used repeatedly in order to rewrite Eq. (30.) in the following form

$$\begin{aligned}
&\begin{bmatrix} H^T S H & K_0 & C_{w0}^T \\ K_0 & 0 & C_{u0}^T \\ C_{w0} & C_{u0} & \gamma R \end{bmatrix} \begin{Bmatrix} \Delta\mathbf{u} \\ \Delta\mathbf{w} \\ \Delta\boldsymbol{\alpha} \end{Bmatrix} = \begin{Bmatrix} 0 & H^T S \\ M & 0 \\ 0 & 0 \end{Bmatrix} \begin{Bmatrix} \Delta\boldsymbol{\delta} \\ \Delta\tilde{\mathbf{u}} \end{Bmatrix} \\
&\quad + \begin{bmatrix} A & 0 & 0 \\ 0 & A & 0 \\ 0 & 0 & \Lambda^T \end{bmatrix} \begin{Bmatrix} \Lambda\Delta\boldsymbol{\alpha} \circ A^T \Delta\mathbf{w} \\ \Lambda\Delta\boldsymbol{\alpha} \circ A^T \Delta\mathbf{u} \\ A^T \Delta\mathbf{w} \circ A^T \Delta\mathbf{u} \end{Bmatrix},
\end{aligned} \qquad (34.)$$

where subscripts 0 denote matrices related to $u_0$, $w_0$ and $\alpha_0$. In particular,

$$C_{u0} = \Lambda^T \mathrm{diag}(A^T u_0) A^T, \quad C_{w0} = \Lambda^T \mathrm{diag}(A^T w_0) A^T. \qquad (35.)$$

Eq. (34) can be written in the compact form

$$K_h \Delta\mathbf{u}_h = M_h \Delta\boldsymbol{\delta}_h + A_h \Theta(A_h^T \Delta\mathbf{u}_h), \qquad (36.)$$

which emphasizes the direct relationship between uncertainties of the given data $\Delta\boldsymbol{\delta}$, $\Delta\tilde{\mathbf{u}}$ and those of the unknown vectors $\Delta\mathbf{u}$, $\Delta\mathbf{w}$, $\Delta\boldsymbol{\alpha}$. Here, $K_h$, $M_h$, $A_h$ are known scalar matrices, and $\Delta\mathbf{u}_h$ depends upon the unknown interval vectors $\Delta\mathbf{u}$, $\Delta\mathbf{w}$ and $\Delta\boldsymbol{\alpha}$. Further, $\Delta\boldsymbol{\delta}_h$ depends upon the known interval vectors $\Delta\boldsymbol{\delta}$ and $\Delta\tilde{\mathbf{u}}$. $A_h^T \Delta\mathbf{u}_h$ is composed of the secondary unknown vectors $A^T\Delta\mathbf{u}$, $A^T\Delta\mathbf{w}$ and $\Lambda\Delta\boldsymbol{\alpha}$. The functional $\Theta(\ )$ in Eq. (34) maps $A_h^T \Delta\mathbf{u}_h$ into the following interval vector

$$\Theta(A_h^T \Delta\mathbf{u}_h) = \Theta\left(\begin{Bmatrix} A^T \Delta\mathbf{u} \\ A^T \Delta\mathbf{w} \\ \Lambda\Delta\boldsymbol{\alpha} \end{Bmatrix}\right) = \begin{Bmatrix} \Lambda\Delta\boldsymbol{\alpha} \circ A^T \Delta\mathbf{w} \\ \Lambda\Delta\boldsymbol{\alpha} \circ A^T \Delta\mathbf{u} \\ A^T \Delta\mathbf{w} \circ A^T \Delta\mathbf{u} \end{Bmatrix}. \qquad (37.)$$

If the square matrix $K_h$ is invertible, Eq. (36.) can be recast into the following fixed-point form

$$\Delta \mathbf{u}_h = \left(K_h^{-1} M_h\right) \Delta \boldsymbol{\delta}_h + \left(K_h^{-1} A_h\right) \Theta\left(A_h^T \Delta \mathbf{u}_h\right) \quad (38.)$$

which is solvable by a new variant of the Neumaier and Pownuk's (2007) method. In particular, we introduce auxiliary variable $\mathbf{v}_h = A_h^T \Delta \mathbf{u}_h$ and the corresponding fixed-point equation follows from Eq. (38) as

$$\mathbf{v}_h = A_h^T \Delta \mathbf{u}_h = \left(A_h^T K_h^{-1} M_h\right) \Delta \boldsymbol{\delta}_h + \left(A_h^T K_h^{-1} A_h\right) \Theta(\mathbf{v}_h) \quad (39.)$$

From this, the following iterative scheme is proposed to find a guaranteed enclosure for $\mathbf{v}_h$. The iteration starts from the trivial initial guess $\mathbf{v}_h^1 = \left(A_h^T K_h^{-1} M_h\right) \Delta \boldsymbol{\delta}_h$ and proceeds in accord with

$$\mathbf{v}_h^{i+1} = \left\{\left(A_h^T K_h^{-1} M_h\right) \Delta \boldsymbol{\delta}_h + \left(A_h^T K_h^{-1} A_h\right) \Theta(\mathbf{v}_h^i)\right\} \cup \mathbf{v}_h^i, \quad (40.)$$

where $\cup$ denotes interval hull of two intervals, and superscripts of $\mathbf{v}_h$ denote iteration steps. The iteration stops when there is no change in $\mathbf{v}_h$ in two consecutive steps, and the converged result is denoted by $\mathbf{v}_h^*$. This is an outer solution for the exact fixed-point $\mathbf{v}_h$ in Eq. (39.), due to the isotonic inclusion of interval operations (Moore, *et al.* 2007). An outer solution for $\Delta \mathbf{u}_h$ is obtained by substituting $A_h^T \Delta \mathbf{u}_h$ in Eq. (38.) with $\mathbf{v}_h^*$. Then the final interval enclosures **u**, **w** and **α** are obtained by subtracting $\Delta \mathbf{u}$, $\Delta \mathbf{w}$ and $\Delta \boldsymbol{\alpha}$ (i.e. $\Delta \mathbf{u}_h$) from $u_0$, $w_0$ and $\alpha_0$ respectively. To further reduce overestimation, the scalar matrices $K_h^{-1} M_h$, $K_h^{-1} A_h$ in Eq. (38.) and $A_h^T K_h^{-1} M_h$, $A_h^T K_h^{-1} A_h$ in Eq. (39.) are calculated before multiplication with the interval vectors $\Delta \boldsymbol{\delta}_h$ and $\Theta(\mathbf{v}_h)$.

## 5 Interval-Based Parameter Identification

In summary, the flowchart of the proposed two-step interval-based inverse algorithm is given in Figure 1. Assume that a finite element model for the structure under study is given. First, we use the deterministic inverse solver introduced in section 3 to estimate a scalar or degenerated interval solution for the unknown parameters. In the second step, the deterministic estimate is used as an initial guess for the interval-based inverse solver defined in section 4. The numerical experiments discussed later on provide strong evidence that the proposed two-step algorithm gives guaranteed interval enclosures of the exact parameters.

Note that the scalar matrices $A$, $\Lambda$ and $M$ are assembled from their element counterparts $A_e$, $\Lambda_e$ and $M_e$, and the constraint matrix $C$ accounts for essential boundary conditions. The interval load uncertainty vector $\boldsymbol{\delta}$ and the measurement vector $\tilde{\mathbf{u}}$ are then determined. Note that $\tilde{\mathbf{u}}$ guarantees to enclose the exact system response, and it is corrupted with random noise. In particular, to simulate realistic conditions, $\tilde{\mathbf{u}}$ is computed as follows:

1. Use a structural FEM model (not necessarily that used in the inversion) to generate the exact measurement data $\tilde{u}_{exact}$.
2. The interval vector $\tilde{\mathbf{u}}_{exact}$ is set with midpoint value $\tilde{u}_{exact}$ and radius equal to the device tolerance $\delta$.
3. An ensemble of perturbed measurements $\tilde{u}_i$ are generated by adding random noise to $\tilde{u}_{exact}$. The random noise is chosen smaller than the tolerance $\delta$ so that $\tilde{u}_i \in \tilde{\mathbf{u}}_{exact}$.
4. Perturbed interval measurement vectors $\tilde{\mathbf{u}}_i$ are generated using $\tilde{u}_i$ as midpoint and device

tolerance $\delta$ as radius. $\tilde{\mathbf{u}}_i$ guarantees to contain $\tilde{u}_{exact}$, i.e. $\tilde{u}_{exact} \in \tilde{\mathbf{u}}_i$.

5. The measurement vector $\tilde{\mathbf{u}}$ is obtained as the intersection of all the $\tilde{\mathbf{u}}_i$ in the ensemble. As a result, $\tilde{\mathbf{u}}$ contains a random perturbation and it still guarantees to contain $\tilde{u}_{exact}$, i.e. $\tilde{u}_{exact} \in \tilde{\mathbf{u}}$.

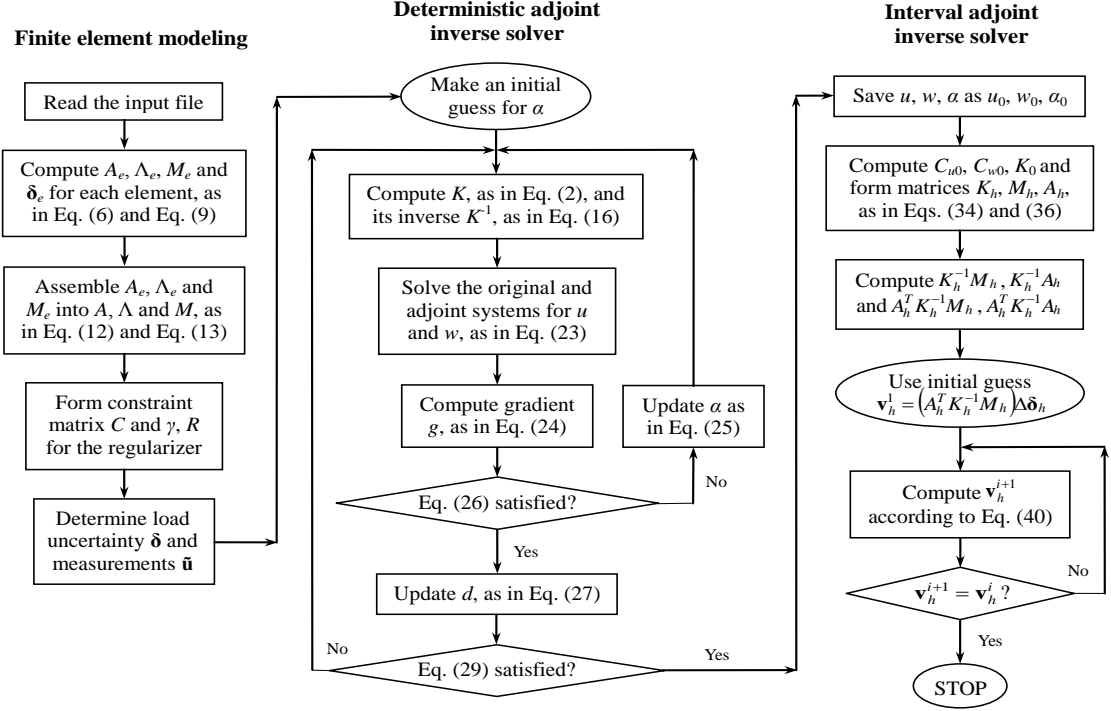

Figure 1. Flowchart for interval-based parameter identification.

In the deterministic solver, to illustrate the robustness of the proposed algorithm, the initial guess is set as $E = 60$ GPa for a structure made of copper, and $E = 160$ GPa for steel for all elements. Then the gradient $g$ in Eq. (24) at the current iteration is computed from the solution vectors $u$ and $w$ of the original and adjoint systems [see Eq. (23)]. Further, we use the weak Wolfe criterion for the inexact line search, setting $\tau_l = 1/4$ and $\tau_u = 1/2$ in Eq. (26.). The Polak-Ribière-Polyak rule in Eq. (28.) is used for the update of the descending directions. In the stopping criterion (29.), the error tolerance $\tau$ is set equal to $1 \times 10^{-10}$ under all circumstances.

In the interval solver, before starting the iteration, we first compute the matrices $C_{u0}$, $C_{w0}$, $K_0$ in Eq. (34). Then we compute the block matrices $K_h$, $M_h$, $A_h$ in Eq. (36), and $K_h^{-1} M_h$, $K_h^{-1} A_h$, $A_h^T K_h^{-1} M_h$, $A_h^T K_h^{-1} A_h$ are computed in advance to solve for $\Delta \mathbf{u}_h$ and $\mathbf{v}_h$ in Eqs. (38) and (39), respectively. As $K_h$, $M_h$, $A_h$ contain a significant number of null-entries, it is more efficient to perform the matrix multiplications and matrix inversions block-by-block. Then the modified version in Eq. (40) of the iterative enclosure method (Newmaier and Pownuk, 2007) is used to compute an enclosure of the unknown parameters $\mathbf{v}_h^*$ from the trivial initial guess $\mathbf{v}_h^1 = \left( A_h^T K_h^{-1} M_h \right) \Delta \boldsymbol{\delta}_h$.

## 6 Numerical Benchmark Problems

The proposed interval inverse algorithm is coded in INTLAB (Rump, 1999), which is an interval arithmetic extension package developed for the MATLAB environment. To test the performance of the method, we consider parameter identification of the Young's moduli of *i)* a fixed-end bar, *ii)* a truss, *iii)* a simply supported beam, and *iv)* a planar frame. Our numerical results show that the proposed method is able to provide a guaranteed interval enclosure of the exact parameters.

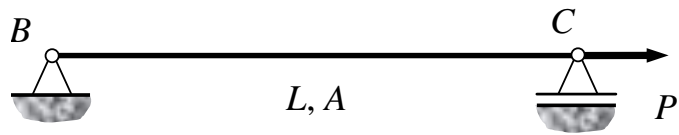

Figure 2. A fixed-end bar subject to concentrated traction at the other end.

## 6.1 Fixed-end bar

Consider a straight bar of length $L = 5$ m, as shown in Figure 2. The pin-roller bar is subject to concentrated force $P = 100$ kN at one end C. The cross section of the bar is uniform, with an area $A = 0.005$ m$^2$. Only axial deformations are allowed, and the bar is modeled by 10 equal-length planar truss elements with uniform material properties. For each element,

$$E = 115 + 10\sin(7x/L) - 5\cos(17x/L) \text{ GPa}, \quad (41.)$$

where $x$ is the coordinate of element centroid, and the values of $E$ are given up to four significant digits. The same 10-element model is also used to generate measurement data. Axial displacements at 10 equally distributed nodes along the bar are collected into the exact measurement vector $\tilde{u}_{exact}$. The interval measurement vector $\tilde{\mathbf{u}}$ is obtained from 3 sets of perturbed measurements $\tilde{u}_i$ with device tolerance $\pm 2 \times 10^{-6}$ m. The results are listed in Table 1. Note that $\tilde{\mathbf{u}}$ contains $\tilde{u}_{exact}$, and uncertainties in $\tilde{\mathbf{u}}$ range from 0.1% to 2%.

Table 1. Exact and perturbed measurement data for the fixed-end bar of Figure 2. The device tolerance is the same for all measurements, $\pm 2 \times 10^{-6}$ m, and 3 sets of perturbed measurements are sampled to define the perturbed data.

| Node # | $\tilde{u}_{exact}$ ($10^{-3}$ m) | $\tilde{\mathbf{u}}$ ($10^{-3}$ m) | | Difference ($10^{-3}$ m) | | Uncertainty (%) | |
|---|---|---|---|---|---|---|---|
| | | Lower Bound | Upper Bound | Lower Bound | Upper Bound | Lower Bound | Upper Bound |
| 1 | 0.09091 | 0.09042 | 0.09226 | -0.00049 | 0.00135 | -0.534 | 1.482 |
| 2 | 0.17281 | 0.17202 | 0.17570 | -0.00079 | 0.00289 | -0.458 | 1.671 |
| 3 | 0.24991 | 0.24789 | 0.25020 | -0.00202 | 0.00029 | -0.809 | 0.118 |
| 4 | 0.33208 | 0.33197 | 0.33265 | -0.00011 | 0.00057 | -0.032 | 0.171 |
| 5 | 0.41980 | 0.41975 | 0.42022 | -0.00005 | 0.00042 | -0.012 | 0.102 |
| 6 | 0.50713 | 0.50554 | 0.50771 | -0.00159 | 0.00058 | -0.315 | 0.114 |
| 7 | 0.59813 | 0.59800 | 0.60031 | -0.00013 | 0.00218 | -0.021 | 0.365 |
| 8 | 0.69694 | 0.69638 | 0.69975 | -0.00056 | 0.00281 | -0.080 | 0.403 |
| 9 | 0.79119 | 0.79014 | 0.79157 | -0.00105 | 0.00038 | -0.133 | 0.048 |
| 10 | 0.87466 | 0.87357 | 0.87555 | -0.00109 | 0.00089 | -0.125 | 0.101 |

This problem has 10 measurements and 10 unknown element Young's moduli $E_i$, and it has an analytical solution. Since the bar is statically determined, axial forces in each element equal to the concentrated traction $P$ at the free end. Then $\mathbf{E}_i$ depends upon the displacements $\mathbf{u}_i$, $\mathbf{u}_{i-1}$ of the neighboring nodes, viz.

$$N = \mathbf{E}_i A \frac{\mathbf{u}_i - \mathbf{u}_{i-1}}{L_e} \quad \Rightarrow \quad \mathbf{E}_i = \frac{NL_e}{A(\mathbf{u}_i - \mathbf{u}_{i-1})}, \tag{42.}$$

where $N = P = 100$ kN is the axial force, $A$ is the cross section area, $L_e = L/10$ is the element length, and $\mathbf{u}_0 = 0$ accounts for the boundary condition at the hinged end.

The problem is well-posed, so no regularization is required. The initial guess $E = 60$ GPa for all the elements. We needed 60 iterations to reach convergence in the deterministic stage, and 12 iterations in the interval stage. The estimated and exact solutions are plotted in Figure 3. Here, the lower and upper bounds of the estimated solution are the dashed lines with triangular markers, and the exact solution is the solid line with rectangular markers. As one can see, the exact values of the Young Moduli are contained by the interval bounds.

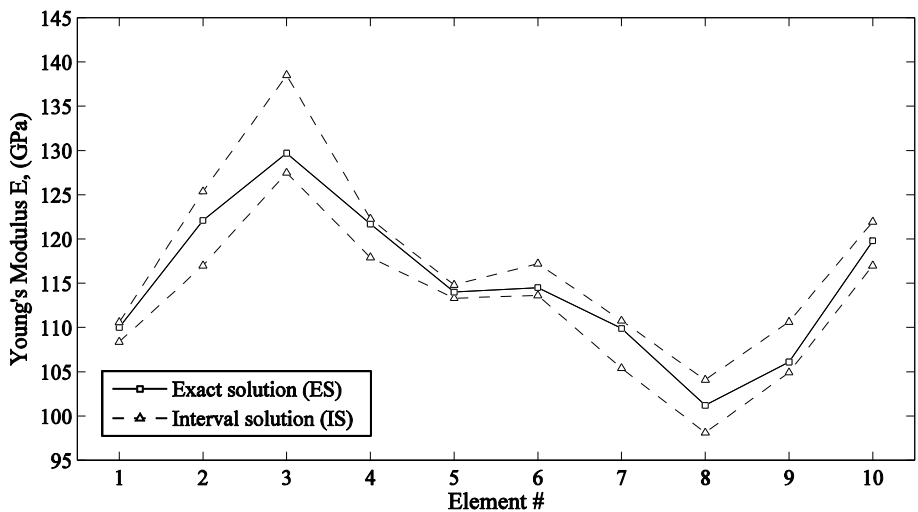

Figure 3. Interval-based identification of Young's moduli of the fixed-end bar of Figure 2: exact values (solid lines with squares) and interval solution (dashed lines with triangles), which is indistinguishable from the Monte Carlo predictions from an ensemble of 10,000 simulations (measurement uncertainty level 0.1-2%).

Table 2 compares the numerical solution $\mathbf{E}_N$ from the proposed method against the analytical solution $\mathbf{E}_A$ from Eq. (42.). The upper bounds of the two solutions are identical, while the lower bounds of $\mathbf{E}_N$ are always smaller than the lower bounds of $\mathbf{E}_A$. In other words, $\mathbf{E}_N$ guarantees to enclose $\mathbf{E}_A$. Exact Young's moduli and relative differences $(E_N - E_A)/E_A \times 100\%$ for the lower and upper bounds of the two interval solutions are also included in the table.

Note that the row of $K_h^{-1} M_h$ corresponding to Young's modulus $\mathbf{E}_i$ of the $i$-th element has all of the entries close to zero, except those at columns corresponding to the measurements $\mathbf{u}_i$ and $\mathbf{u}_{i-1}$ at the neighboring nodes. In addition, the two entries have similar magnitude and opposite sign. This is in agreement with the analytical solution given in Eq. (42.), that is: the modulus $\mathbf{E}_i$ of the $i$-th element is only a function of $\mathbf{u}_i$ and $\mathbf{u}_{i-1}$.

Table 2. Exact Young's moduli and predicted values for the fixed-end bar of Figure 2. Relative differences ($E_N - E_A)/E_A \times 100\%$ for the lower and upper bounds of the two interval solutions are also listed.

| Element # | Exact (GPa) | $\mathbf{E}_N$ (GPa) | | $\mathbf{E}_A$ (GPa) | | Relative Diff. (%) | |
|---|---|---|---|---|---|---|---|
| | | Lower Bound | Upper Bound | Lower Bound | Upper Bound | Lower Bound | Upper Bound |
| 1 | 110.0 | 108.37 | 110.59 | 108.39 | 110.59 | -0.020 | 0.000 |
| 2 | 122.1 | 117.00 | 125.37 | 117.27 | 125.37 | -0.231 | 0.000 |
| 3 | 129.7 | 127.48 | 138.52 | 127.90 | 138.52 | -0.332 | 0.000 |
| 4 | 121.7 | 117.91 | 122.30 | 117.98 | 122.30 | -0.066 | 0.000 |
| 5 | 114.0 | 113.30 | 114.80 | 113.31 | 114.80 | -0.009 | 0.000 |
| 6 | 114.5 | 113.63 | 117.22 | 113.68 | 117.22 | -0.048 | 0.000 |
| 7 | 109.9 | 105.39 | 110.75 | 105.52 | 110.75 | -0.120 | 0.000 |
| 8 | 101.2 | 98.12 | 104.08 | 98.28 | 104.08 | -0.169 | 0.000 |
| 9 | 106.1 | 104.91 | 110.63 | 105.05 | 110.63 | -0.137 | 0.000 |
| 10 | 119.8 | 116.99 | 121.95 | 117.09 | 121.95 | -0.085 | 0.000 |

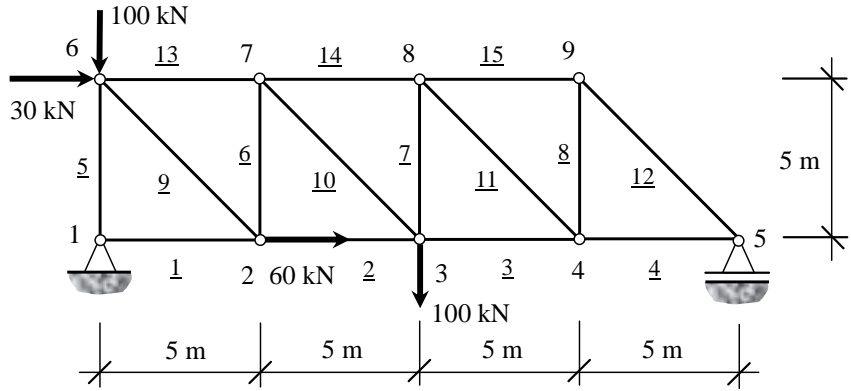

Figure 4. A simply-supported truss subject to concentrated loads.

### 6.2 Simply-supported truss

The second example is a simply supported truss composed of 15 bars, subject to concentrated loads, as shown in Figure 4. Nodes of the truss are labeled from 1 to 9, and the bars are labeled from 1 to 15. Horizontal load 60 kN is applied at node 2, vertical load 100 kN at node 3, horizontal load 30 kN and vertical load 100 kN at node 6. The bars have uniform cross sections with area $A = 0.005$ m$^2$. Each bar is modeled by one planar truss element with constant material property, and the corresponding Young's modulus is listed in the second column of Table 4. Here we assume that bar 3 and 13 are damaged, and their effective Young's moduli are 80 GPa and 60 GPa, respectively.

The same finite element model is used to compute the exact measurement data. To illustrate the performance of the current method under different forms of measurements, nodal displacements of bottom nodes 2 to 5, as well as strains of medium-height bars 5 to 12, are measured. The device tolerance is $\pm 1 \times 10^{-5}$ m for nodal displacement measurements, and $\pm 1 \times 10^{-6}$ for strain measurements. The measurement vector $\tilde{\mathbf{u}}$ is obtained from 3 sets of perturbed $\tilde{u}_i$, and the results are shown in Table 3.

The uncertainties in **ũ** range from 0.06% to 2%, approximately.

Table 3. Exact and perturbed measurements for the simply-supported truss of Figure 4. The device tolerance is $\pm 1\times 10^{-5}$ m for nodal displacements, and $\pm 1\times 10^{-6}$ for strains. 3 sets of perturbed measurements are sampled to yield the perturbed data.

|  | Exact ($10^{-3}$ m) | Perturbed data ($10^{-3}$ m) | | | | | Exact ($10^{-4}$) | Perturbed data ($10^{-4}$) | | | |
|---|---|---|---|---|---|---|---|---|---|---|---|
|  |  | Lower Bound | Uncertainty (%) | Upper Bound | Uncertainty (%) |  |  | Lower Bound | Uncertainty (%) | Upper Bound | Uncertainty (%) |
| $u_2$ | 0.7557 | 0.7532 | -0.321 | 0.7586 | 0.382 | $\varepsilon_5$ | -2.3246 | -2.3306 | -0.256 | -2.3210 | 0.155 |
| $v_2$ | -5.1714 | -5.1732 | -0.036 | -5.1591 | 0.238 | $\varepsilon_6$ | -0.6822 | -0.6827 | -0.078 | -0.6674 | 2.161 |
| $u_3$ | 1.3922 | 1.3871 | -0.369 | 1.4021 | 0.711 | $\varepsilon_7$ | 0.9664 | 0.9661 | -0.025 | 0.9777 | 1.167 |
| $v_3$ | -7.6368 | -7.6393 | -0.032 | -7.6349 | 0.025 | $\varepsilon_8$ | 1.0388 | 1.0309 | -0.769 | 1.0456 | 0.648 |
| $u_4$ | 2.8297 | 2.8141 | -0.551 | 2.8310 | 0.047 | $\varepsilon_9$ | 1.1427 | 1.1387 | -0.346 | 1.1457 | 0.265 |
| $v_4$ | -4.3003 | -4.3045 | -0.097 | -4.2914 | 0.208 | $\varepsilon_{10}$ | 1.1028 | 1.1000 | -0.253 | 1.1043 | 0.133 |
| $u_5$ | 3.2930 | 3.2924 | -0.019 | 3.3089 | 0.482 | $\varepsilon_{11}$ | -1.4241 | -1.4354 | -0.795 | -1.4213 | 0.199 |
|  |  |  |  |  |  | $\varepsilon_{12}$ | -1.3736 | -1.3748 | -0.088 | -1.3591 | 1.058 |

This problem has 15 measurement and 15 unknowns. It is well-posed and no regularizer is needed. The initial guess $E = 60$ GPa is used. 465 iterations are run in the deterministic stage, and 12 iterations in the interval stage. In Table 4 and Figure 5, the obtained interval solution (IS) is compared against the exact solution (ES) and Monte Carlo (MC) predictions based on an ensemble of 10,000 simulations. In each simulation $k$, a random measurement vector $\tilde{u}_k$ is chosen within the interval bounds of $\tilde{\mathbf{u}}$, i.e. $\tilde{u}_k \in \tilde{\mathbf{u}}$. The corresponding solution $\alpha_k$ is obtained from the deterministic inverse solver formulated in section 3, and the Monte Carlo solution $\boldsymbol{\alpha}_{MC}$ is given by the minimum and maximum values of all $\alpha_k$ in the ensemble, that is $\boldsymbol{\alpha}_{MC} = [\min_k \alpha_k, \max_k \alpha_k]$.

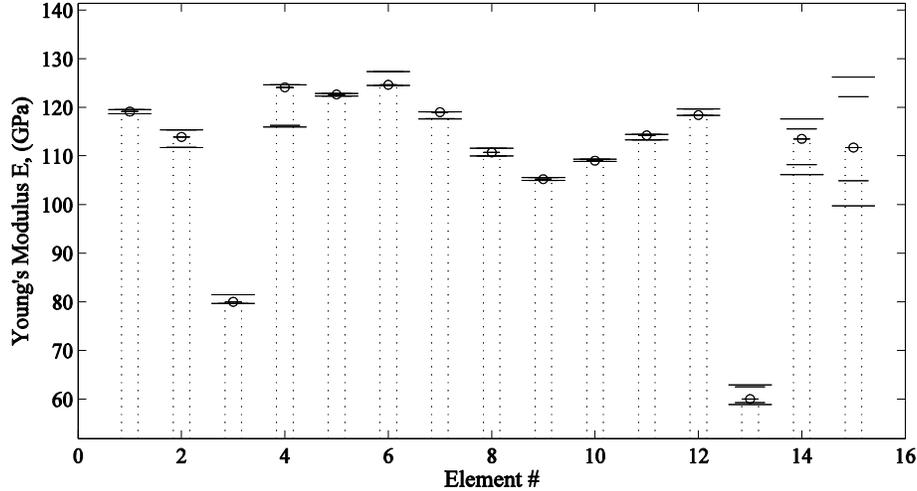

Figure 5. Interval-based identification of Young's moduli of a simply-supported truss of Figure 4: short bars with circular markers denote the exact values; the long bars denote interval prediction from the proposed method; median-length bars with circles denote Monte Carlo predictions from an ensemble of 10,000 simulations.

Clearly, both IS and MC predictions enclose the exact values of the Young's moduli, and IS contains MC. It is observed that the interval enclosures of IS are very tight for elements 5 to 12, and very wide for elements 13 to 15. This can be explained as follows: The structure in Figure 4 is statically determined, and the estimate of the Young's modulus $\mathbf{E}_i$ ($i = 1,\ldots,15$) is directly related to the element axial strain $\varepsilon_i$. In particular, for the elements 5 to 12, the axial strains $\varepsilon_i$ ($i = 5,\ldots,12$) are directly measured. Thus the corresponding uncertainty level in Young's moduli $\mathbf{E}_i$ ($i = 5,\ldots,12$) is relatively small. For elements 1 to 4, the axial strains $\varepsilon_i$ ($i = 1,\ldots,4$) are indirectly obtained from measured nodal displacements of the neighboring nodes, $\varepsilon_i = (\mathbf{u}_{i+1} - \mathbf{u}_i)/L_e$, where $L_e$ is the element length and $\mathbf{u}_1 = 0$ accounts for the boundary condition. Because of the interval subtraction, the uncertainty level in $\mathbf{E}_i$ ($i = 1,\ldots,4$) is larger. For elements 13 to 15, the axial strains $\varepsilon_i$ ($i = 13, 14, 15$) are indirectly related to multiple components of the measurement data. As an example, consider element 13. By checking entries of $K_h^{-1}M_h$ in Eq. (38), $\mathbf{E}_{13}$ is related to the nodal displacements $\mathbf{u}_2$, $\mathbf{v}_2$, $\mathbf{u}_3$, and $\mathbf{v}_2$, as well as the axial strains $\varepsilon_5$, $\varepsilon_6$, $\varepsilon_9$, $\varepsilon_{10}$. As a result, the uncertainty level in $\mathbf{E}_i$ ($i = 13, 14, 15$) is the largest.

Table 4. Exact and predicted Young's modulus for the simply-supported truss of Figure 4. Relative error of the interval solutions from the proposed method and Monte Carlo predictions from an ensemble of 10,000 simulations.

| Element # | Exact (GPa) | Proposed method (GPa) | | | | Monte Carlo method (GPa) | | | |
|---|---|---|---|---|---|---|---|---|---|
| | | Lower Bound | Uncertainty (%) | Upper Bound | Uncertainty (%) | Lower Bound | Uncertainty (%) | Upper Bound | Uncertainty (%) |
| 1 | 119.1 | 118.64 | -0.383 | 119.48 | 0.322 | 118.65 | -0.381 | 119.48 | 0.322 |
| 2 | 113.9 | 111.68 | -1.949 | 115.35 | 1.277 | 111.75 | -1.884 | 115.33 | 1.257 |
| 3 | 80.0 | 79.62 | -0.474 | 81.44 | 1.804 | 79.66 | -0.431 | 81.42 | 1.780 |
| 4 | 124.1 | 115.92 | -6.593 | 124.63 | 0.425 | 116.24 | -6.333 | 124.59 | 0.392 |
| 5 | 122.6 | 122.29 | -0.256 | 122.79 | 0.156 | 122.29 | -0.255 | 122.79 | 0.156 |
| 6 | 124.6 | 124.47 | -0.105 | 127.35 | 2.210 | 124.50 | -0.077 | 127.35 | 2.208 |
| 7 | 119.0 | 117.62 | -1.163 | 119.03 | 0.027 | 117.63 | -1.153 | 119.03 | 0.025 |
| 8 | 110.7 | 109.97 | -0.655 | 111.56 | 0.776 | 109.99 | -0.642 | 111.56 | 0.775 |
| 9 | 105.2 | 104.92 | -0.267 | 105.57 | 0.348 | 104.92 | -0.265 | 105.57 | 0.347 |
| 10 | 109.0 | 108.85 | -0.136 | 109.28 | 0.255 | 108.85 | -0.133 | 109.28 | 0.253 |
| 11 | 114.2 | 113.29 | -0.796 | 114.43 | 0.202 | 113.30 | -0.788 | 114.43 | 0.199 |
| 12 | 118.4 | 118.29 | -0.096 | 119.67 | 1.071 | 118.30 | -0.087 | 119.67 | 1.069 |
| 13 | 60.0 | 58.84 | -1.931 | 62.89 | 4.818 | 59.24 | -1.272 | 62.51 | 4.178 |
| 14 | 113.5 | 106.16 | -6.466 | 117.56 | 3.578 | 108.15 | -4.710 | 115.57 | 1.823 |
| 15 | 111.7 | 99.72 | -10.726 | 126.23 | 13.011 | 104.88 | -6.107 | 122.14 | 9.349 |

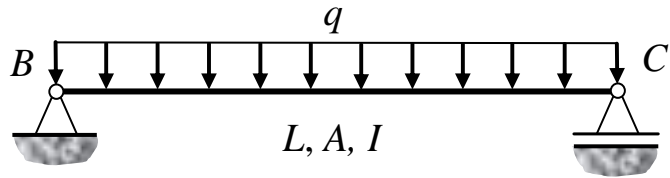

Figure 6. A simply-supported beam subject to uniformly distributed vertical load.

*6.3 Simply-supported beam*

The third example is a simply-supported beam subject to uniformly distributed vertical load $q = 100$ kN/m, as shown in Figure 6. The beam has a length $L = 2$ m, and a 5 cm×3 cm rectangular cross section (cross section area $A = 0.015$ m$^2$ and moment of inertia $I = 1.125 \times 10^{-4}$ m$^4$). The beam is subject to lateral deformation, and 20 two-node Euler-Bernoulli beam elements are used in the finite element mesh. The stiffness matrix is computed using the three-node Gaussian quadrature rule. In order to generated a continuous material field, Young's moduli at the quadrature nodes are linearly interpolated from those at the nodal values

$$E = 220 + 10\sin(6x/L) - 5\cos(13x/L) \text{ GPa}, \qquad (43.)$$

where $x$ is the nodal coordinate, and the values are given up to four significant digits. The parameter vector $\boldsymbol{\alpha}$ has 21 components, one for each mesh node.

In the first case, a finer 80-element finite element model is used to generate the measurement data. Young's moduli are linearly interpolated from the abovementioned 21-node material mesh. Further, 9 lateral deflections at equidistant points along the beam are collected as measurements. The measurement vector $\tilde{\mathbf{u}}$, which has 9 components, is obtained from 3 sets of perturbed data $\tilde{u}_i$ with device tolerance $\pm 2 \times 10^{-6}$ m. The resulting $\tilde{\mathbf{u}}$ has uncertainties ranging from 0.1% to 1%, and contains the exact measurement data.

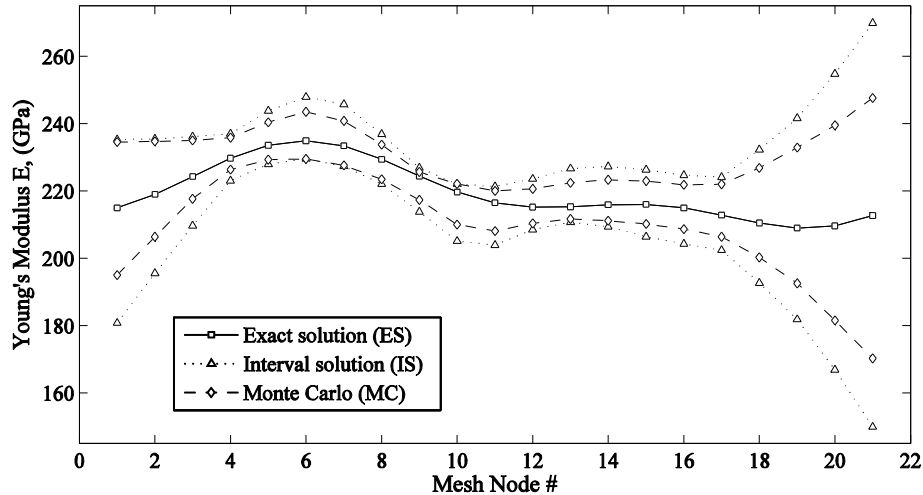

Figure 7. Interval-based identification of the Young's moduli of the simply-supported beam of Figure 6 under uniformly distributed load: interval solution (IS), exact solution (ES) and Monte Carlo (MC) prediction from an ensemble of 10,000 simulations (measurement uncertainty level 0.1-1%).

The problem is ill-posed, since only 9 measurements are available to estimate 21 unknown parameters. This requires regularization. The regularizer weight $\gamma$ should be chosen with caution: it has to be large enough to avoid useless estimate or even divergence with unbounded intervals, but not that large, otherwise the solution will be over-smoothen (Hansen, 2010). Here, we use a second-order regularization matrix $R$ and $\gamma = 1 \times 10^{-3}$. For the proposed method, the initial guess $E = 160$ GPa for all elements. Convergence is attained in 289 and 37 iterations in the deterministic and interval stages, respectively. The interval estimates are compared against the exact Young's modulus from Eq. (43.) and Monte Carlo prediction from an ensemble of 10,000 simulations. Figure 7 shows the exact solution (ES, solid lines with rectangular markers), the interval solution (IS, dotted lines with triangular markers) and the Monte Carlo prediction (MC, dashed lines with diamond markers). Observe that IS indicates a high level of uncertainty near both ends, especially near the right end, which is attributed to the relatively small bending moment near the ends. In addition, both IS and MC guarantees to enclose ES everywhere, and IS contains MC.

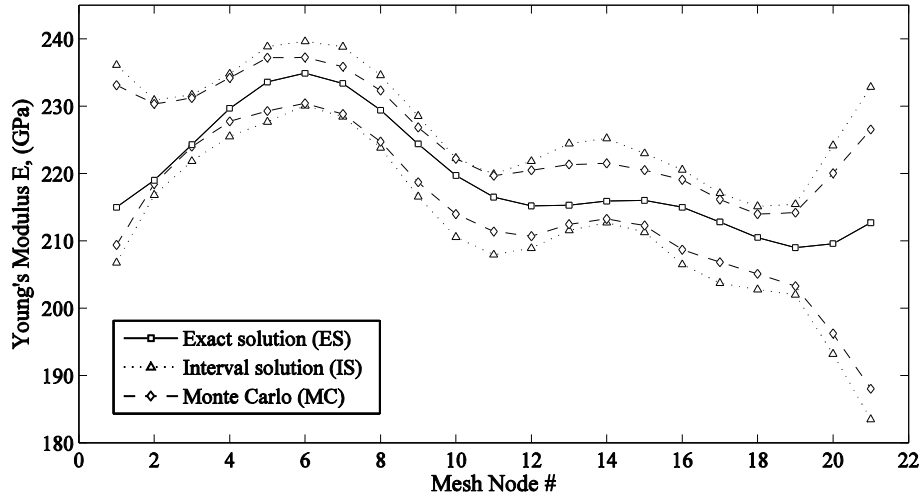

Figure 8. Interval-based identification of Young's moduli of the simply-supported beam of Figure 6 under uniformly distributed load and bending moments at both ends: interval solution (IS), exact solution (ES) and Monte Carlo (MC) prediction from an ensemble of 10,000 simulations (measurement uncertainty level 0.1-1%).

In the second case, two opposing bending moments $M = 50$ kN·m are added to the ends B and C, in order to create a more uniform bending moment diagram for the beam. In addition, rotation angles $\theta_B$ and $\theta_C$ at both ends are measured. The device tolerance is now $\pm 5 \times 10^{-6}$ m for deflections and $\pm 2 \times 10^{-5}$ rad for $\theta_B$ and $\theta_C$. As a result, the level of uncertainty in $\tilde{\mathbf{u}}$ ranges from 0.1% to 1%, roughly the same as in the first case. IS and MC predictions are compared against the exact values ES in Figure 8. Note that the level of uncertainty at the ends is reduced significantly. This is due to the additional bending moments at the ends and extra measurements $\theta_B$ and $\theta_C$. Indeed, the maximum level of uncertainty at the ends is approximately 13% on the left and 23% on the right. In the previous case of Figure 7 the uncertainty levels increase to approximately 25% on the left and 56% on the right. Near the mid-span, the level of uncertainty is slightly reduced from about 8% in Figure 7 to about 5% in Figure 8.

Finally, we point out that interval solutions guarantee to enclose all possible predictions associated with different probabilistic distributions of the measurements, either symmetrical or not (see Figure 9).

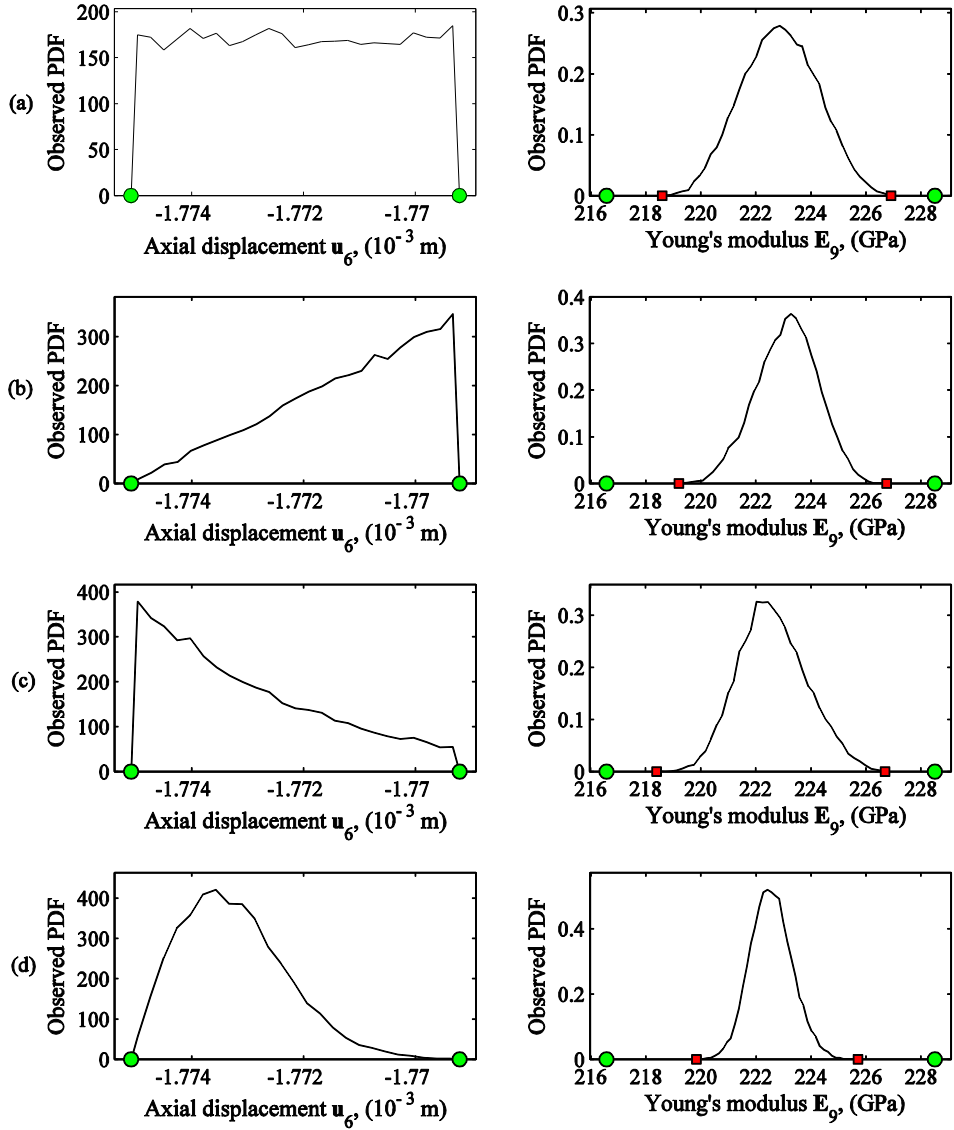

Figure 9. Comparison between the interval solution and Monte Carlo prediction of the Young's modulus $\mathbf{E}_9$ of the simply-supported beam of Figure 6 from an ensemble of 10,000 simulations: (left) observed probability density function (PDF) of axial displacement measurements $\mathbf{u}_6$ sampled from (a) uniform, (b) triangular, (c) truncated exponential and (d) truncated Rayleigh probability distributions (interval endpoints denoted by circular markers); (right) corresponding observed PDF of the Young modulus $\mathbf{E}_9$, interval solution (endpoints denoted by circular markers) and Monte Carlo predicted interval [min($E_9$) max($E_9$)] (square markers).

*6.4 Two-bay two-story frame*

The fourth example is a two-bay two-story planar frame hinged to the ground, subject to uniformly distributed vertical loads on each floor, as shown in Figure 9. The frame is composed of six columns and four beams, labeled as $C_j$ ($j = 1,…,6$) and $B_j$ ($j = 1,…,4$), respectively. Connecting joints and supports are labeled as nodes 1 to 9. Uniformly distributed vertical loads $q_i$ ($i = 1,…,4$) are applied on $B_i$, where $q_1 = q_2 = 109.45$ kN/m and $q_3 = q_4 = 51.08$ kN/m.

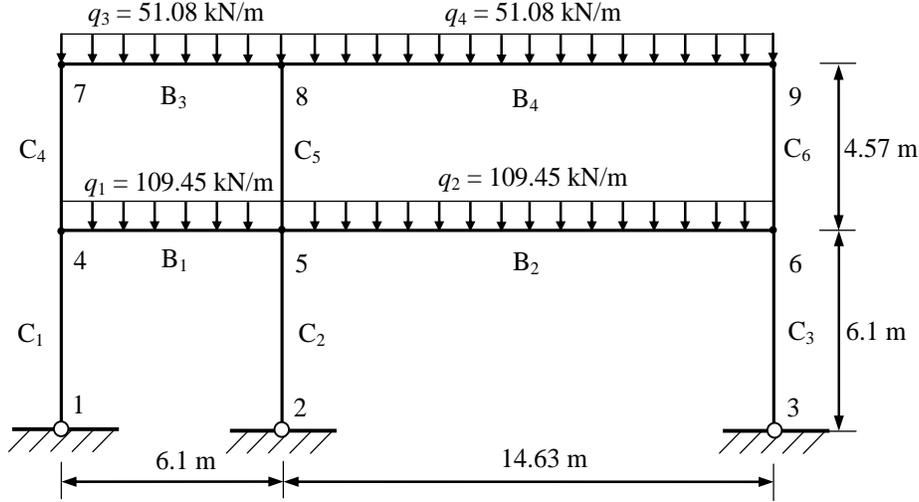

Figure 9. A two-bay two-story frame subject to uniformly-distributed vertical loads on each floors.

Each member of the frame has uniform cross section and material property. The corresponding cross section area $A$, moment of inertia $I$ and Young's modulus $E$ are listed in Table 5. Ten two-node Euler-Bernoulli beam elements are used to model the frame, one for each member.

Measurement data used in the inverse algorithm is generated from the same 10-element finite element model. Only nodal displacement $u_i$, $v_i$ and rotation angle $\theta_i$ at nodes 4 to 9 ($i = 4,\ldots,9$) are included in the measurement vector $\tilde{\mathbf{u}}$. $\tilde{\mathbf{u}}$ is obtained from 3 sets of perturbed measurements $\tilde{u}_i$, and the corresponding device tolerance is $\pm 2\times 10^{-5}$ m for nodal displacements and $\pm 2\times 10^{-5}$ rad for rotation angles. The level of uncertainty in $\tilde{\mathbf{u}}$ ranges from approximately 0.1% to 1%, with the exception of $\theta_4 = [-1.2442, -0.9825]\times 10^{-4}$ rad (22.2% uncertainty).

Table 5. Geometric and material properties for the members of the two-bay two-story frame shown in Figure 10.

|   | Shape | $A$ ($10^{-4}$ m$^2$) | $I$ ($10^{-8}$ m$^4$) | $E$ (GPa) |   | Shape | $A$ ($10^{-4}$ m$^2$) | $I$ ($10^{-8}$ m$^4$) | $E$ (GPa) |
|---|---|---|---|---|---|---|---|---|---|
| $C_1$ | W12×19 | 35.940 | 5411.00 | 210 | $B_1$ | W27×84 | 160.000 | 118625.96 | 205 |
| $C_2$ | W14×132 | 250.320 | 63683.41 | 214 | $B_2$ | W36×135 | 256.130 | 324660.51 | 208 |
| $C_3$ | W14×109 | 206.450 | 51612.70 | 205 | $B_3$ | W18×40 | 76.130 | 25473.36 | 215 |
| $C_4$ | W10×12 | 22.835 | 2239.32 | 201 | $B_4$ | W27×94 | 178.710 | 136107.68 | 214 |
| $C_5$ | W14×109 | 206.450 | 51612.70 | 204 |   |   |   |   |   |
| $C_6$ | W14×109 | 206.450 | 51612.70 | 206 |   |   |   |   |   |

In this benchmark case, 18 measurements (6 nodes × 3 DOF) are used to predict the Young's modulus $E$ of the 10 members. The problem is well-posed and no regularizer is required. Initial guess $E = 160$ GPa is used. The results are compared with the exact solution and the Monte Carlo solution with 10,000 runs in Figure 10, following the same guidelines as in Figure 5 of the simply-supported truss. Observe that the interval solution provides a guaranteed enclosure of both the exact and Monte Carlo solutions.

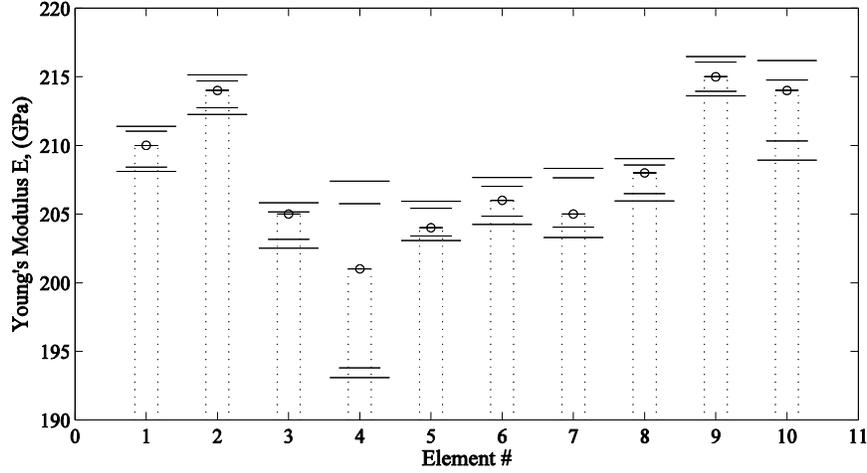

Figure 10. Interval-based identification of Young's moduli of the two-bay two-story frame in Figure 10 subject to uniformly distributed loads: short bars with circular markers denote the exact values; long bars denote interval predictions from the proposed method; median-length bars denote the Monte Carlo predictions from an ensemble of 10,000 simulations (measurement uncertainty level 0.1%-1%).

In Figure 10, note that the width of the interval estimate $\mathbf{E}_4$ for the Young's modulus of the left column $C_4$ on the upper floor, is much wider than other estimates. By examining the entries of $K_h^{-1}M_h$ in Eq. (38.), we note that the wide enclosure is mainly caused by the lateral displacements $\mathbf{v}_4$ and $\mathbf{v}_7$ at nodes 4 and 7. These two measurements are modeled by two intervals with about 1% uncertainty, i.e. $\mathbf{v}_4$ = [–2.3599, –2.3399]×10$^{-3}$ m and $\mathbf{v}_7$ = [–3.4548, –3.4186]×10$^{-3}$ m. In order to obtain a narrower interval prediction for $\mathbf{E}_4$, we increase the accuracy of the measurements $\mathbf{v}_4$ and $\mathbf{v}_7$, and reduce the level of uncertainty to about 0.2%, i.e. $\mathbf{v}_4$ = [–2.3515, –2.3465]×10$^{-3}$ m and $\mathbf{v}_7$ = [–3.4378, –3.4288]×10$^{-3}$ m. The results are depicted in **Error! Reference source not found.**, showing a significant increase in the accuracy of the predicted value for $\mathbf{E}_4$. In particular, the previous estimate in Figure 10 is $\mathbf{E}_4$ = [193.09, 207.39] GPa (7.1% uncertainty), and that in **Error! Reference source not found.** is $\mathbf{E}_4$ = [197.72, 203.34] GPa (2.8% uncertainty).

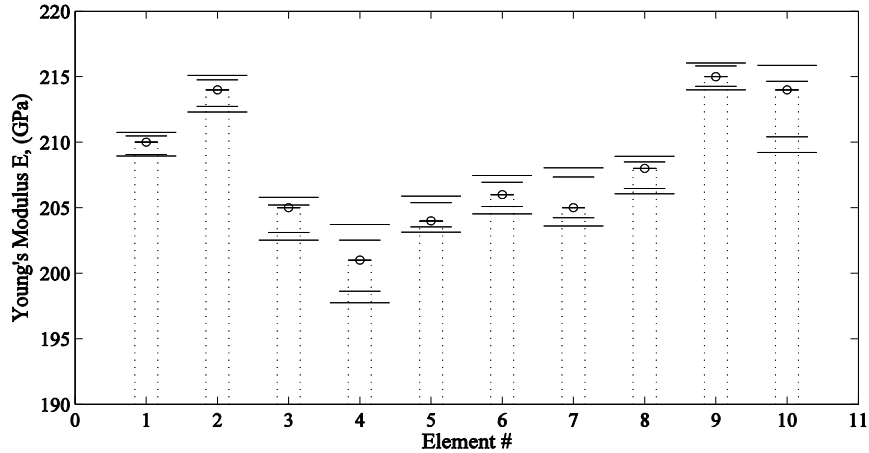

Figure 11. Interval-based identification of the Young's moduli of the two-bay two-story frame in Figure 10 subject to uniformly distributed loads using more accurate measurements in $\mathbf{v}_4$ and $\mathbf{v}_7$ (uncertainty level 0.2%) than those used to obtain the estimates shown in Figure 10: short bars with circular markers denote the exact values; long bars denote interval prediction from the proposed method; median-length bars denote the Monte Carlo prediction from an ensemble of 10,000 simulations.

# 7   Conclusions

An interval-based parameter identification is proposed for structural static problems. Uncertainties in the system are modeled by intervals, and IFEM is exploited to handle uncertainties. The proposed inverse algorithm stems from an adjoint-based optimization formulation, and it provides an interval estimate of the unknown parameters (e.g. Young's moduli). The associated nonlinear interval equations are solved by means of a new variant of the iterative enclosure method. In addition, overestimation is reduced by means of a new decomposition of the IFEM matrices **K** and **f**, which limits multiple occurrences of the same variable in the IFEM equations by separating deterministic and interval terms. The interval solution from the proposed solver guarantees to enclose the exact parameters, as confirmed by several numerical benchmark problems, and it always contains Monte Carlo predictions.